\documentclass[12pt]{article}

\usepackage{amsmath}

\usepackage{amssymb}

\usepackage{graphicx}

\def \UN{\hbox{ \rm 1\hskip -3.2pt I}}

\def\P{\mathbb{P}}
\def\E{\mathbb{E}} %

\def\R{\mathbb{R}}

\def\S{{\mathcal S^{n-1}}} 
\def\Var{\mathop{\rm Var}}

\newcommand{\findem}{\hfill\hbox{\hskip 4pt
\vrule width 5pt height 6pt depth 1.5pt}\vspace{.5cm}\par}

\newtheorem{theo}{Theorem}

 \newtheorem{lem}{Lemma}

\title{The tail of the maximum of smooth Gaussian fields on fractal sets }

\author{Jean-Marc Aza\"{\i}s
\thanks{%
Universit\'e de Toulouse, IMT ESP,
Universit\'e Paul Sabatier,  31062 Toulouse Cedex 9. France, jean-marc.azais@math.univ-toulouse.fr }
\and Mario Wschebor \thanks{%
Centro de Matem\'atica. Facultad de Ciencias. Universidad de la
Rep\'ublica. Calle Igua 4225. 11400 Montevideo. Uruguay,
wschebor@cmat.edu.uy}}

\begin{document}

 \maketitle

\emph{Short Title: } Maximum on fractals.


\begin{abstract}

We study the probability distribution of the maximum $M_S $ of a smooth stationary Gaussian field defined on a fractal subset $S$ of $\R^n$. Our main result is the equivalent of the asymptotic behavior of the tail of the distribution $\P(M_S>u)$ as $u\rightarrow +\infty.$ The basic tool is Rice formula for the moments of the number of local maxima of a random field.
\end{abstract}

\section{Introduction}

Let  $\{ X(t) , t \in \R \}$ be a stationary centered Gaussian process with constant variance equal to $1$ and $ \mathcal{C}^1$ paths  and let $M$ be its maximum on a compact interval, for example $[0,1]$.
It is well-known since the pioneering work of S.O. Rice \cite{rice} that
$$
\P\{M>u\}  \leq  1-\Phi(u)   + (const)~\varphi (u),
$$
 where $\varphi $ is the density of the standard Gaussian law,  $\Phi$ its cumulative distribution function and $ (const) $ is a constant depending upon the law of the process. In what follows, $(const)$ will denote a constant that may change in each occurrence.\\

 Under stronger hypotheses, one can prove the equivalence (see, for example,  Piterbarg  \cite{piter}, Aza\"{\i}s and Wschebor \cite{AW09}  Prop 4.1 and 4.2).

$$
\P\{M>u\}  \simeq (const)~\varphi (u) ~~~\text{as}~~~u\rightarrow +\infty.
$$
We are using  the following definitions and notation:
 For $f,g$ real-valued functions of one real variable, we put ($x_0$ can be a real number, $+\infty$ or $-\infty $):
$$
f\simeq g~~\text{as}~~x\rightarrow x_0
$$
if $\lim_{x\rightarrow x_0} f(x)/g(x)=1.$\\

Similarly, we put
$$
f\sim g~~\text{as}~~x\rightarrow x_0
$$
if there exists a positive constant $c$ such that $1/c \leq f(x)/g(x) \leq c$ for all $x \neq x_0$ in some neighborhood of $x_0.$\\

Let now $\mathcal{X}=\{X(t), t \in \R^n \}$, $n\geq 1$  be a sufficiently smooth real-valued centered stationary Gaussian random field, with a law that satisfies certain non-degeneracy assumptions. We consider a compact set $S,~S\subset \R^n$ having topological dimension $d$ and possessing a geometry that satisfies a certain number of smoothness properties that we do not detail here. Let $M_S$ be the maximum of the random field on $S$. It is known (See for example Adler \cite{A81}, Piterbarg \cite{Pit96}, Taylor, Takemura and Adler \cite{TTA}, Adler and Taylor \cite{AT07} Aza\"{\i}s and Wschebor \cite{AW08}, \cite{AW09}) that:\\

\begin{equation}\label{f:d}
\P\{M_S>u\}   \simeq   (const)~u^{d-1} \varphi (u)~~~\text{as}~~~u\rightarrow +\infty.
\end{equation}

\bigskip

Under certain conditions, it is possible to go beyond this equivalence and obtain a whole expansion for $u\rightarrow +\infty$. In fact, the approximation becomes very accurate, in the sense that for large values of $u$, the error term is exponentially smaller than the principal term. Other extensions are known, for example to non-stationary Gaussian fields (\cite{AT07}, \cite{AW09}, \cite{Pit96}).\\

The goal of this paper is to extend formula (\ref{f:d}) assuming that $S$ is a fractal set, in which case no results have been published yet, as far as the authors know. We will assume in the sequel that $S$ is a bounded subset of $\R^n$ and denote $M_S=\sup \{X(t):t\in S\}$. Assuming that the paths $t\rightsquigarrow X(t)$ are continuous, it is easy to see that $M_S$ is actually a random variable, without any further hypothesis on $S$.\\

We will prove that whenever $S$ is Minkowski measurable in the strong sense (Falconer \cite{Fal90}, definition below) formula (\ref{f:d}) holds true replacing $d$ by the Minkowski dimension and the constant factor can be explicitly computed from the law of the process and the content of the domain. This includes a large class of self-similar sets. In other cases, we are able to obtain the asymptotic behavior of $\P (M_S>u)$ as $u\rightarrow +\infty$ up to a constant factor or - more weakly - the equivalent of $~\log \P (M_S>u)$.\\

Another interesting situation is when $S$ which has an outer Minkowski content but may have a non-smooth boundary (see example 4 below). In this case, we obtain the asymptotic behavior of $\P (M_S>u)$ for large $u$, which depends upon the Lebesgue measure and the perimeter of the set $S$.\\

We are only considering here equivalences of the form (\ref{f:d}) not full expansions which appear to be essentially more complicated when the set $S$ is fractal. This will be the object of future work.\\

For each subset $S$ of $\R^n$ and each $\epsilon >0$, we will denote $S^{\epsilon}$  the tube of radius $\epsilon $ around $S$, that is the set $S^{\epsilon}=\{t\in \R^n:~dist(t,S)\leq \epsilon \}$, where $dist$ denotes Euclidean distance.  This set is  sometimes called "parallel set"  or "$\epsilon$-neigbourhood". The set $S$ is called ``Minkowski measurable in the strong sense'', with Minkowski dimension $d$ and content $C$, if

$$
\lambda _n( S^ {\epsilon})  \simeq C \epsilon ^{n-d} ~~ \text{as}~~ \epsilon \rightarrow 0,
$$

\noindent $\lambda _n$ denotes Lebesgue measure in $\R^n.$\\

Similarly, $S$ is called ``Minkowski measurable in the weak sense'', with Minkowski dimension $d$ if
$$
\lambda _n( S^ {\epsilon})  \sim  \epsilon ^{n-d} ~~ \text{as}~~ \epsilon \rightarrow 0.
$$
More generally, the Minkowski dimension $d$ is defined as:
$$
d=n-\lim_{\epsilon \rightarrow 0}\frac{\log \lambda_n(S^{\epsilon})}{\log\epsilon},
$$
whenever the limit exists.\\

\section{Main theorem}\label{statproof}

Let $\tilde{ B } $ be a ball centered at the origin such that $S\subset \tilde{B}$. We set $ B = \tilde{ B }^{1}$.\\

\noindent For $u\in \R$ consider the excursion set
$$
E_u =\{t \in B:X(t) \geq u\}.
$$
Let $t$ be a local maximum of $X(.)$, that is,
$$
X(t) \geq X(s)~~\forall~s\in B\cap U~~\text{for some open ball}~U~\text{centered at }t.
$$
If $t\in E_u$, we denote by $K_t $ the connected component  of  $E_u$ that  contains $t$. Note that such a $t$ can lie on the interior $ \dot{B}$  of the ball $B$ or on its boundary $\partial B$.\\


We make the  following assumptions on the random field $\mathcal{X}:$

\begin{itemize}
\item  It is Gaussian, centered and stationary. We denote its covariance by \\$\Gamma(t)=\E (X(s)X(s+t))$.  
\item  The paths  are of class $\mathcal{C}^3$.
\item For  all $s\neq t \in B$ the random vector:
 $$
( X(s), X(t),X'(s),X'(t)) ~~\mbox{has a non degenerate distribution}.
 $$
 \item For all  $t\in B,~\lambda \in \S$ the random vector:
 $$
( X(0),X'(0),X''(0)\lambda) ~~\mbox{has a non degenerate distribution}.
 $$
\item  With no loss of generality we assume in addition that it is centered, and for every $t$ one has $\Var(X(t)) =1$, $\Var(X'(t)) =I_n$, where $I_n$ is the identity matrix. Concerning this last condition, in any case, one can make a linear transformation of the parameter space so that it is satisfied, changing the set $S$ into its image under this transformation.
\end{itemize}

  If  $\lambda$ is a unit vector in $\S$, we denote $X''_{\lambda} (s)=\lambda^\top X''(s)\lambda $. $B(t,r)~(t\in \R^n,~r>0 )$ denotes the open ball with center $t$ and radius $r$ in $\R^n.$ \\

The following lemma  is essential for the proof of the main theorem.

 \begin{lem}
Let $\alpha$  be a real number $0 < \alpha<1$.  Let  $M_u $  be the  number of local maxima above $u$ in $\dot B$. Under our assumptions, the following  quantities are ``negligible''  meaning by that that they are $o( u^{-1}\varphi(u))$ as $u\rightarrow +\infty$:
  \begin{enumerate}
\item  $\P (A_1)$ where $A_1=\{ \exists  \mbox{ a local maximum  in } B  \mbox { with value }\geq u + 1\} $.
\item   $\E( M_u(M_u-1)) $.
\item  $\P (A_2) $ where $A_2=\{ \exists ~ \mbox{ two or more local maxima in } \dot{B} \mbox { with value }\geq u  \} $.
\item  $\P (A_3)$ where
\begin{equation*}
\aligned
A_3&=\{ \exists ~ \mbox{a local maximum } t\in \dot{B}\\
&\mbox { such that } u<X(t)<u+1,~\min \{X''_{\lambda} (s): s \in B(t,\rho), \lambda = \frac{s-t }{\|s-t\|} \}\leq - X(t) - u^\alpha  \}
\endaligned
\end{equation*}
where
$$
\rho=\rho (u)=u^{-\beta} \ \
\mbox{  with  } \beta >(1-\alpha )/2
$$
This means that the connected component associated to $t$ is not too small.

\item  $\P (A_4)$ where
\begin{equation*}
\aligned
A_4&=\{ \exists ~ \mbox{a local maximum } t\in B\\
&\mbox { such that } u<X(t)<u+1,~\max \{X''_{\lambda} (s):s \in B(t,\rho), \lambda = \frac{s-t }{\|s-t\|} \} \geq - X(t) + u^\alpha  \}
\endaligned
\end{equation*}
This means that the connected component associated to $t$ is not too large.
\end{enumerate}

\end{lem}

\emph{Proof : }

\noindent 1 is a direct consequence of the fact  that $E(M_u) \leq (const) u^{n-1} \varphi(u)$ \\

\noindent 2 is proved in \cite{Piterbarg} and in  Delmas' PHD dissertation \cite{delmas} under slightly more restrictive  conditions on the random field, but the  proof goes through with the same arguments under our hypotheses.  The result is also quoted in \cite{AD}. \\

\noindent  3 is a direct consequence of 2. \\

\noindent The proof of 4 is very similar but simpler than that of 5.\\

\noindent \emph{Proof of 5}. Let $  M_4$ be the number of local maxima in the interior of $  B$  such that\\ $u<X(t) <u+1$  and
$$
\max \{X''_{\lambda} (s) : s \in B(t,\rho), \lambda = \frac{s-t }{\|s-t\|} \}\geq - X(t) + u^\alpha
$$
Then $\P(A_4) \leq \E(M_4) $. This expectation can be computed by means of a Rice formula (see \cite{AW09}, Chap. 6):
\begin{equation*}
\aligned
&\E(M_4)\\
 = &\int_u^{u+1}dx
\int_B    \E \Big( |\det[X''(t)]|\UN_{X''(t)\prec 0}\UN_{\max\{X''_{\lambda} (s): s \in B(t,\rho), \lambda = \frac{s-t }{\|s-t\|}\}  \geq - x  +u^\alpha} \big| X(t)=x,X'(t)=0\Big)\\
&p_{X(t),X'(t)}(x,0)\lambda_n(dt).
\endaligned
\end{equation*}

\noindent $\UN_A$ denotes the indicator function of the set $A$, $\prec$ means negative definite,  $ p_{\zeta}(.) $ is the probability density function of the random vector $\zeta$ .\\

Using the stationarity of the law of the random field and the form of the density, we get:
\begin{equation}\label{esp3}
\aligned
&\E(M_4)\\
&\leq (const)\int_u^{u+1}\varphi (x)
 \E \Big( |\det[X''(0)]|\UN_{\max\{X''_{\lambda} (s): s \in B(0,\rho), \lambda = \frac{s}{\|s\|} \} \geq - x +u^\alpha  } \big| X(0)=x,X'(0)=0\Big)dx
\endaligned
\end{equation}
We use the  Schwarz inequality for the conditional expectation:
\begin{equation}\label{cauchy}
\aligned
&\Bigg[\E \Big( |\det[X''(0)]|\UN_{\max\{X''_{\lambda} (s): s \in B(0,\rho), \lambda = \frac{s}{\|s\|}\} \geq - x + u^\alpha  } \big| X(0)=x,X'(0)=0\Big)\Bigg]^2\\
&\leq \E \Big( |\det[X''(0)]|^2  \big| X(0)=x,X'(0)=0\Big)\\
&~~~~~~~~~~~~~~~~~\times \P \Big(\max\{X''_{\lambda} (s): s \in B(0,\rho), \lambda = \frac{s }{\|s\|}  \}\geq - x +u^\alpha \big| X(0)=x,X'(0)=0 \Big).
\endaligned
\end{equation}
For the first factor in the right-hand side of (\ref{cauchy}), note that $X(0)$ and $X'(0)$ are independent as well as  $X''(0) $ and $ X'(0)$. Regression formulas imply that, under the condition,

$$
X''(0)=-x I_n+\zeta,
$$
where the matrix $\zeta=X''(0)+X(0)I_n$ is independent of $X(0),X'(0).$ Hence:
$$
\E \Big( |\det[X''(0)]|^2  \big| X(0)=x,X'(0)=0\Big) \leq (const)\big[x^{2n}+1 \big].
$$
Let us turn to the second factor in the right-hand side of (\ref{cauchy}).
Under the condition, we have:
$$
X''_{\lambda} (s)=\Gamma''_\lambda (s) x+ R_\lambda (s),
$$
where $R_\lambda (s)$  is a random field with parameter $s$ and bounded variance and  $\Gamma$ denotes the covariance  function of  $ \mathcal{X}$ which is of class $ \mathcal{C}^6$. Using standard properties of the covariance and the continuity of its fourth derivative, we can write
$$
\Gamma''_\lambda (s) \leq -1 + (const)~\|s\|^2.
$$
Taking into account our choice for $\alpha$ and $\beta$, it follows that for large enough $u$ we get the following bound for the second factor:
\begin{equation}\label{proba}
\aligned
&\P \Big( \max \{ R_\lambda (s) : s \in B(0,\rho) \}\geq   u^\alpha   - (const) \|s\|^2 (u+1) \Big)\\
&\leq
\P \Big( \min \{ R_\lambda (s):s \in B(0,\rho) \}\leq  - (const) u^\alpha  \Big)
\\
&\leq (const)~\exp \big[ -(const)u^{2\alpha}\big].
\endaligned
\end{equation}

The last inequality follows from the classical Landau-Shepp-Marcus-Fernique inequality for the tails of the supremum of a.s. bounded Gaussian sequences (see for example \cite{AW09},(2.33)). Replacing in (\ref{esp3}) we finish the proof of 5. \findem

\begin{theo}\label{main}
Assume that the process $\mathcal{X} $ satisfies the above conditions. Then,
\begin{enumerate}
\item If the parameter set $S$ is Minkowski measurable in the strong sense, with dimension $d$ and content $C$, one has, as $u\rightarrow +\infty$:
\begin{equation}\label{minkstrong}
\P (M_S>u)\simeq \frac{C}{2^{d/2}\pi ^{n/2}}\Gamma \big( 1+(n-d)/2 \big)~u^{d-1}\varphi (u).
\end{equation}
\item If the parameter set $S$ is Minkowski measurable in the weak sense, with dimension $d$, one has, as $u\rightarrow +\infty$:
\begin{equation}\label{minkweak}
\P (M_S>u)\sim u^{d-1}\varphi (u).
\end{equation}
\item If the Minkowski dimension $d$ of the parameter set $S$ is well-defined then, as $u\rightarrow +\infty$:
\begin{equation}\label{minksolo}
\log \P (M_S>u) =-\frac{1}{2}u^2+(d-1+o(1))\log u.
\end{equation}
\end{enumerate}
\end{theo}

\emph{Proof :} We denote
$$
A=\bigcup_{i=1}^4A_i
$$

According to the Lemma, as $u\rightarrow +\infty $
$$
\P(A)=o(\varphi (u)/u).
$$
If $u$ is large enough, whenever  the complement set $A^C$ occurs and the excursion set is non-empty, the connected component $K_t$  associated to the local maximum  $t \in  B$  satisfies
\begin{equation}\label{f:kitu}
 B\big(t, \underline{r}\big)
   \subset  K_t \subset
   B\big(t, \bar{r}\big)
 \end{equation}
Where
$$
\underline{r}=\sqrt{2\frac{ X(t) -u}{ X(t)+u^\alpha}  },~~~\bar{r}=\sqrt{2\frac{X(t)-u}{ u-u^\alpha} }.
$$

 This means that for $u$ sufficiently large,  such that $\sqrt{ 2/(u-u^\alpha)} <1$ the component associated to a local maximum  belonging to  $\partial B$    does not touch $S$. This implies that  maxima in $\partial B$ do not have to be taken into account when considering  the event\\$\{M_S>u\} \cap A^C$.

Using Lemma 1,  on $A^C$  the excursion set $E_u$ has at most one connected component associated to a  local maximum  in $\dot{B}$.\\

Let us prove (\ref{f:kitu}). Notice first that if $u$ is large enough, one has $B(t,\underline{r}),B(t,\bar{r})\subset B(t,\rho)$. We prove the first inclusion in (\ref{f:kitu}), the proof of the second one is similar.\\

Let $s\in B(t,\underline{r}), s\neq t.$ Write the Taylor expansion
$$
X(s)=X(t)+\frac{1}{2}\|s-t\|^2~X''_{\lambda}(\theta)
$$
where $\lambda = (s-t)/\|s-t\|$ and $\theta$ is some point in the linear segment $[s,t]$. So,
$$
X(s)\geq X(t)-\frac{1}{2}\|s-t\|^2~(X(t)+u^{\alpha})\geq X(t)-\frac{X(t)-u}{X(t)+u^{\alpha}}(X(t)+u^{\alpha})=u.
$$
This proves that $s\in K_t.$\\

Now consider:
\begin{equation}\label{chain}
\aligned
\P(M_S>u) = o(\varphi (u)/u)+\P \Big(\{\exists t\in B, X(.) &~\text{has a local maximum at}~t, u<X(t)<u+1 \}\\
&\cap \{K_t\cap S \neq \emptyset \}\cap A^C \Big).
\endaligned
\end{equation}
%


Then, an instant reflection shows that we have the following upper and lower bounds for $\P (M_S>u)$:
\begin{equation}\label{uplo}
\aligned
&\P(M_S>u) \leq o(\varphi (u)/u)+ \P \big(\{\exists t\in\dot{ B,} X(.)~\text{has a local maximum at}~t,
u<X(t)<u+1,t\in S^{\bar{r}} \}\big),\\
&\P(M_S>u) \geq o(\varphi (u)/u)+ \P \big(\{\exists t\in\dot{ B}, X(.)~\text{has a local maximum at}~t,
u<X(t)<u+1,t\in S^{\underline{r}} \}\big).
\endaligned
\end{equation}
Denote
$$
M_{u,u+1,r}=\#\{t \in\dot{ B} :t\in S^r,~ X(.) ~\text{has a local maximum at}~ t, u<X(t)<u+1\}.
$$
The first inequality in (\ref{uplo}) implies:
$$
\P(M_S>u) \leq o(\varphi (u)/u)+ \E \big( M_{u,u+1,\bar{r}} \big).
$$
As for the second inequality in (\ref{uplo}), we have:
$$
0 \leq \E \big( M_{u,u+1,\underline{r}}~ \big)-\P \big(\{M_{u,u+1,\underline{r}}\geq 1 \big)
\leq \frac{1}{2}\E \big( M_{u,u+1} (M_{u,u+1}-1)\big)=o(\varphi (u)/u)
$$
so that
$$
\P(M_S>u) \geq o(\varphi (u)/u)+\E \big( M_{u,u+1,\underline{r}} \big)
$$


\medskip

To finish the proof we study the behavior as $u\rightarrow +\infty$ of  $\E (M_{u,u+1,\bar{r}}) $ and $\E (M_{u,u+1,\underline{r}}) $. Notice that $\underline{r},\bar{r}$ are random.
Let us consider case 1., when the set $S$ is Minkowski measurable. Cases 2 and 3 are entirely similar.\\

Let us look at $\E (M_{u,u+1,\bar{r}}) $. Our tool is again Rice formula (See \cite{AW09}, Chap. 6) which in this case takes the form:

\begin{align}\label{ricemark}
&\E (M_{u,u+1,\bar{r}})   \notag \\
&=\int_u^{u+1}dx\int_B \E \Big( |\det[X''(t)]|\UN _{\{X''(t)\prec 0\}}\UN_{\{t\in S^{\bar{r}}\}} \Big/ X(t)=x,X'(t)=0\Big)p_{X(t),X'(t)}(x,0)\lambda_n(dt) \notag \\
 & =I
 \end{align}

We have, for $u$  large enough:

\begin{align}
I&=\int_u^{u+1}dx\int_B \UN_{\{t\in S^{\bar{r}^*}\}}\E \Big( |\det[X''(t)]|\UN _{\{X''(t)\prec 0\}} \Big|X(t)=x,X'(t)=0\Big) \notag \\
&~~~~~~~~~~~~~~~~~~~~~~~~~~~~~~~~~~~~~~~~~~~~~~~~~~~~~~~~~~~~~~~~~~~~~~~~~~~~~~~~~~~~~\frac{1}{(2\pi)^{n/2}}
\varphi (x)\lambda_n(dt) \notag \\
&=\frac{1}{(2\pi)^{n/2}}\int_u^{u+1} \lambda_n \big(S^{\bar{r}^*}\big)\E \Big( |\det[X''(0)]|\UN _{\{X''(0)\prec 0\}} \Big| X(0)=x,X'(0)=0\Big)\varphi (x)dx \label{rice4}
\end{align}
where we have replaced the density by its value and used that the conditional expectation is independent of $t$, since the process is stationary. $\bar{r}^*$ is defined as the value of $\bar{r}$ when replacing $X(t)$ by $x$. Note that it is non-random.\\

It is well known that the conditional expectation in the right-hand side of (\ref{rice4}) is equivalent to $x^n$ as $x\rightarrow +\infty.$ (See for example \cite{A81}). Use now that the set $S$ is Minkowski measurable in the strong sense, with dimension $d$ and content $C$. This easily implies that as $u\rightarrow +\infty$ one has:

\begin{equation*}
\aligned
I &\simeq \frac{1}{(2\pi)^{n/2}} \frac{1}{\sqrt{2\pi}}\int_u^{u+1}~x^n\exp (-x^2/2) C
\bigg[ \frac{2(x-u)}{u-u^{\alpha}}\bigg]^{(n-d)/2}dx\\
&\simeq \frac{C}{2^{d/2}\pi ^{n/2}}\frac{1}{\sqrt{2\pi}}~u^{(n+d)/2}\int_u^{u+1}~(x-u)^{(n-d)/2}\exp (-x^2/2)dx
\endaligned
\end{equation*}
Now perform the change of variables $x=u(1+y/u^2)$ in the last integral. We get:

\begin{align*}
I\ & \simeq  \frac{C}{2^{d/2}\pi ^{n/2}}u^{d-1}\varphi (u)\int _0^u y^{(n-d)/2}\exp \big[-y-\frac{y^2}{2u^2}\big] dy\\
&\simeq \frac{C}{2^{d/2}\pi ^{n/2}}u^{d-1}\varphi (u)\Gamma \big(\frac{n-d}{2}+1\big),
\end{align*}
where the last equivalence follows using Lebesgue's Theorem.\\

For the lower bound of $ \P(M_S>u) $, we use again Rice formula to express $\E (M_{u,u+1,\underline{r}}) $ as an integral. The computation is entirely similar and we obtain the same equivalent as $u\rightarrow +\infty$. This ends the proof. \findem

\bigskip

\section{Examples}\label{exam}

\begin{enumerate}
\item \textbf{Smooth manifolds} \\
Our first example is the well-known case in which $S$ is a compact smooth manifold without boundary having topological dimension $d$ embedded in $\R^n.$ In this case, $S$ is Minkowski measurable in the strong sense, with dimension $d$ and content equal to
$$
\sigma_d(S) Vol(B_{n-d}),
$$
where $\sigma_d(S)$ denotes the geometric measure of $S$ and $Vol(B_{n-d}) $ is the volume of the unit ball in $\R^{n-d}$. Replacing in (\ref{minkstrong}) we get:
$$
\P (M_S>u)\simeq \sigma_d(S)\frac{1}{(2\pi)^{d/2}}~u^{d-1}\varphi (u)
$$

This is a much weaker version of results in \cite{TTA} or \cite{AW08}

\item \textbf{Self-similar sets}.

Let $f_i:\R^n\rightarrow \R^n ~~(i=1,..,k)$ be $k$ similarities of $\R^n$, that is,
$$
\|f_i(x)-f_i(y) \| =c_i\|x-y\|~~\forall~~x.y \in \R^n
$$
and we assume that $0<c_i <1$ for $i=1,...,k,~k\geq 2.$ Let $S$ be the self-similar set associated to $f_1,...,f_k$, which means that
$$
S=\bigcup_{i=1}^kf_i(S).
$$
We also assume that the $f_i$'s satisfy the ``open set condition'', that is there exists a non-empty bounded open set $V$ such that
$$
\bigcup_{i=1}^kf_i(V)\subset V
$$
where the union is disjoint.\\

Under these conditions, one can prove that the Hausdorff and Minkowski dimensions of $S$ are equal, and their common value $d$ is the unique solution of the Moran equation:
$$
\sum_{i=1}^kc_i^d=1.
$$

So, if $S$ is self-similar and satisfies the open set condition, one has the logarithmic expansion (\ref{minksolo}).

Classical examples of self-similar sets satisfying these conditions are the homogeneous Cantor sets, the Sierpinski gasket, the von Koch and modified von Koch curves. For a proof of the above statement and a number of other examples and extensions, see \cite{Fal90} ).

\item \textbf{Special self-similar sets}.\\

For self-similar sets one can obtain Minkowski measurability in the strong sense, at the cost of adding a certain number of conditions. One can compute the dimension and the content and then obtain the equivalence for the tail of the distribution of the maximum having the form (\ref{minkstrong}) in the statement of the main theorem . We restrict ourselves to dimension $n=1$. The formulas become more complicated for higher dimensions, but it is possible to extend the results, adding a certain number of conditions (see \cite{DKO10}).\\

In dimension $1$, the result is given in \cite{LP93} and \cite{Fal95}.
We use the notation of
the previous example and assume the same hypotheses. Let the interval $I$ be the convex hull of $S$. Then, for $i\neq i'$, $f_i(I)$ and $f_{i'}(I)$ are disjoint, except possibly at the endpoints, so that $I\setminus \cup _{i=1}^k f_i(I)$ is the union of $k',~k'\leq k-1$ disjoint intervals with lengths $\ell_1,...,\ell_{k'} $.\\

The set $S$ is said to be of lattice-type if the additive subgroup of the reals $\sum_{i=1}^k(\log c_i)\mathbb{Z}$ is discrete. Then, the result is that $S$ is Minkowski-measurable in the strong sense, if and only if it is not of lattice-type, and in that case, the content is given by:
$$
C=\frac{2^{1-d}\sum_{j=1}^{k'}\ell _j^d}{d(1-d)\sum_{i=1}^k c_i^d \log(1/c_i)}.
$$
As an illustrative example, assume that $k=2.$ Then, $S$ is not of lattice type if and only if $\log c_1 /\log c_2$ is irrational.

\end{enumerate}

%

%

\section{Sets with a finite perimeter}

Let $S\subset \R^n$ be a closed set. Whenever the limit
$$
\lim_{\epsilon \rightarrow 0}\frac{\lambda_n(S^{\epsilon} \backslash S)}{\epsilon}=OM(S)
$$
exists, $S$ is said to have the ``outer Minkowski content $OM(S)$.''\\

For a review of this subject and related ones, see for example \cite{ACV}. Let us mention here some general examples of classes of sets having outer Minkowski content:
\begin{enumerate}
\item $S$ is compact with Lipschitz boundary, in which case
\begin{equation}\label{lipboun}
OM(S)=DG(S)=H^{n-1}(\partial S)
\end{equation}
where $DG $ denotes De Giorgi's perimeter (see \cite{DG}) and $H^{\alpha}$ the usual Hausdorff measure of exponent $\alpha$.
\item Let $S\subset \R^n$ be a closed set and define:
$$
U(S)=\{a\in \R^n: \exists ~\text{a unique}~ x\in S ~\text{such that}~dist (a,S)=dist(a,x)\}
$$
$$
reach(S)=\inf_{a\in S}\sup\{r>0:B(a,r)\subset U(S)\}.
$$
$S$ is said to have ``positive reach'' if $reach(S)>0$ and in this case, $S$ has an outer Minkowski content (see \cite{Fed}). Compact convex sets with non-empty interior have positive reach and in this case, (\ref{lipboun}) is also verified, allowing to compute the outer Minkowski content as a perimeter of the Hausdorff measure of the boundary.
\item The outer Minkowski content need not coincide with the Hausdorff measure of the boundary. More precisely, it is possible to construct compact sets $S\subset \R^n$ with non-empty interior, positive reach, hence having outer Minkowski content, but such that $OM(S)>H^{n-1}(\partial S)$ (see for example \cite{ACV}, example 1, where this is done for $n=2$).
\end{enumerate}

We have the following new result:

\begin{theo}

Assume that $S$ has outer Minkowski content. Then, as $u\rightarrow +\infty$:
\begin{equation*}
\P (M_S>u) =  \lambda_n (S)\frac{1}{(2\pi)^{n/2}}~u^{n-1}\varphi (u) + \frac{ OM(S) }{2^{(n+1)/2} \pi ^{(n-1)/2}} u^{n-2} \varphi(u)\big(1+o(1)\big)
\end{equation*}
\end{theo}

The proof follows the one of Theorem \ref{main}, with the only change that one needs a somewhat more precise result on the asymptotic behavior of the conditional expectation. In fact, one can use that (see for example \cite{delmas})
$$
\E \Big( |\det[X''(0)]|\UN _{\{X''(0)\prec 0\}} \Big| X(0)=x,X'(0)=0\Big) = x^n + O(x^{n-2})
$$
and this suffices to adapt the proof of the main theorem to this example.

 \end{document}